\documentclass{amsproc}

\newtheorem{theorem}{Theorem}[section]

\theoremstyle{definition}
\newtheorem{definition}[theorem]{Definition}
\newtheorem{example}[theorem]{Example}

\theoremstyle{remark}

\numberwithin{equation}{section}

\def\rsd{\mbox{$\mathop{\mathrel\times\joinrel\mathrel{\vrule height 5pt
 depth 0pt}} $}}

\def\lsd{\mbox{$\mathop{\mathrel{\vrule height 5pt
 depth 0pt}\joinrel\mathrel\times} $}}

\begin{document}

\title{Examples of masas in C*-algebras}
\author{Jean Renault}
\address{D\'epartment de Math\'ematiques, Universit\'e d'Orl\'eans,
BP 6759, 45067 Orl\'eans Cedex 2, France}
\email{Jean.Renault@univ-orleans.fr}


\date{November 30, 2008 and, in revised form, ?? ??, 2008.}



\begin{abstract}
This paper illustrates the notion of a Cartan subalgebra in a C*-algebra through a number of examples and counterexamples. Some of these examples have a geometrical flavour and are related to orbifolds and non-Hausdorff manifolds.
\end{abstract}

\maketitle
\markboth{Jean Renault}
{Examples}


.

\section{Cartan subalgebras of C$^*$-algebras}

I recall in this section some definitions and results of \cite{ren:cartan}, to which I refer the reader for a more complete exposition.

\begin{definition} We shall say that an abelian
sub-C$^*$-algebra $B$ of a  C$^*$-algebra $A$
is a {\it Cartan subalgebra} if
\begin{enumerate}
\item $B$ contains an approximate unit of $A$;
\item $B$ is maximal abelian;
\item $B$ is regular;
\item there exists a faithful conditional expectation $P$ of $A$ onto $B$.
\end{enumerate}
\end{definition}

The main result of \cite{ren:cartan} is a C$^*$-algebraic version of Feldman-Moore's well known theorem \cite{fm:relations II} on Cartan subalgebras in von Neumann algebras. This theorem establishes an equivalence of categories between twisted countable standard measured equivalence relations and Cartan subalgebras in von Neumann algebras on separable Hilbert spaces. A notable difference in the topological case is that equivalence relations (also known as principal groupoids) have to be replaced by topologically principal groupoids. The definition of a topologically principal groupoid is related to the definition of a topologically free action, as given in \cite[Definition 2.1]{tom:dynamics}.

\begin{definition}\label{principal} We say that a  groupoid $G$ on a topological space $X$ (this means that $G$ has $X=G^{(0)}$ as its unit space) is
{\it topologically principal} if the set of points of
$X$ with trivial isotropy is dense.
\end{definition}

All our examples of topologically principal groupoids will be groupoids of germs. Suppose that an inverse semi-group $\Gamma$ acts on a topological space $X$ by partial homeomorphisms (i.e. homeomorphisms of an open subset of $X$ onto another open subset). Then the set $G$ of germs of these homeomorphisms form a groupoid: let us write $[g(x),g,x]$ the germ of $g\in\Gamma$ at a point $x$ in the domain of $g$; the groupoid operations are simply
$$[x,g,y][y,h,z]=[x,gh,z]\qquad [x,g,y]^{-1}=[y,g^{-1},x].$$
We identify $X$ with $G^{(0)}$ through the map $x\mapsto [x,id,x]$, where $id$ is the identity map of $X$. 
The topology of germs turn $G$ into an \'etale topological groupoid (\'etale means that the range and source maps are local homeomorphisms). Such a groupoid is not necessarily Hausdorff, even if the unit space $X$ is Hausdorff, as we usually assume. If $\Gamma$ is countable and $G^{(0)}$ is a Baire space, then $G$ is topologically principal. Conversely, let $G$ be an \'etale topological groupoid on a topological space $X$. Then, the inverse semi-group of its open bisections acts on $X$ by partial homeomorphisms. The corresponding groupoid of germs is a quotient of $G$. One says that $G$ is {\it effective} if this quotient map is injective. If $G$ is Hausdorff and topologically principal, then it is effective.

Let me also recall the construction of the reduced C$^*$-algebra of a Hausdorff  locally compact groupoid $G$ equipped with a left Haar system $\lambda=\{\lambda^x\}$. When $G$ is \'etale, one uses the counting measures on the fibers $G^x=r^{-1}(x)$ as a Haar system. The following operations turn the space $C_c(G)$ of compactly supported complex-valued continuous functions on $G$ into an involutive algebra:
$$f*g(\gamma)=\int f(\gamma\gamma')g(\gamma'^{-1})d\lambda^{s(\gamma)}(\gamma');$$
$$f^*(\gamma)=\overline{f(\gamma^{-1})}.$$
For each $x\in G^{(0)}$, one defines the representation $\pi_x$ of $C_c(G)$ on the Hilbert space $L^2(G_x,\lambda_x)$, where $G_x=s^{-1}(x)$ and $\lambda_x=(\lambda^x)^{-1}$, by $\pi_x(f)\xi=f*\xi$. One defines the reduced norm $\|f\|_r=\sup\|\pi_x(f)\|$. The reduced C*-algebra $C_r^*(G)$ is the completion of $C_c(G)$ for the reduced norm. We shall need a slight generalization of the above construction. A twist over  a groupoid $G$ is a groupoid extension
$${\bf T}\times X\rightarrow E\rightarrow G$$
where ${\bf T}$ is the circle group, $X$ is a space and, at the level of the unit spaces, the maps
$X\rightarrow E^{(0)}\rightarrow G^{(0)}$  
are identification maps. In the topological setting, we require the maps to be continuous and the identification maps to be homeomorphisms.
We replace the complex-valued functions by the sections of the associated complex line bundle. Essentially the same formulas as above provide the C*-algebra $C_r^*(G,E)$.

\begin{theorem}\label{Cartan}\cite{ren:cartan} Let $(G,E)$ be a twist with $G$ \'etale, second countable locally compact Hausdorff and topologically principal. Then $C_0(G^{(0)})$ is a Cartan subalgebra of $C_r^*(G,E)$.

Conversely, let $B$ be a Cartan sub-algebra of a separable C*-algebra $A$. Then, there exists a twist $(G,E)$ with $G$ \'etale, second countable locally compact Hausdorff and topologically principal and an isomorphism of  $C_r^*(G,E)$ onto $A$ carrying $C_0(G^{(0)})$ onto $B$.
\end{theorem}

\medskip
This theorem extends a theorem of Kumjian \cite{kum:diagonals} who deals with the principal case and introduces the stronger notion of a diagonal. It requires the property of unique extension of states, which is studied in \cite{gre:thesis, abg:extensions}.

\begin{definition}
One says that a sub-C*-algebra $B$ of a C*-algebra $A$ has {\it the unique extension property} if all pure states of $B$ extend uniquely to pure states of $A$. A Cartan subalgebra which has the unique extension property is called a {\it diagonal} (or a {\it diagonal subalgebra}).
\end{definition}

\begin{theorem} \cite{kum:diagonals, ren:cartan} Let $B$ be a Cartan sub-algebra of a separable C*-algebra $A$. Let $(G,E)$ be the associated twist. Then,
$G$ is principal if and only if $B$ has the unique extension property.
\end{theorem}

\section{Examples}

\subsection{Two non-conjugate Cartan subalgebras}

Here are examples of C$^*$-algebras containing at least two non-conjugate diagonal subalgebras (which are not even isomorphic as algebras). I owe the first one to A. Kumjian.

1. Let ${\bf T}$ be the circle and let $n$ be an integer not smaller than 2. Let $D_n$ be the subalgebra of diagonal matrices in the algebra of matrices $M_n({\bf C})$. The C$^*$-algebra $C({\bf T})\otimes M_n({\bf C})$ obviously contains $C({\bf T})\otimes D_n=C({\bf T})^n$ as a diagonal subalgebra. However, as shown for example in \cite [Example 3(iii)]{kum:prelim}, $C({\bf T})\otimes M_n({\bf C})$ can be realized as the crossed product C$^*$-algebra $C({\bf T})\rsd\,\, {\bf Z}_n$, of the action of ${\bf Z}_n={\bf Z}/n{\bf Z}$ on the circle ${\bf T}$ by the rotation of angle $2\pi/n$. Therefore, it also contains $C({\bf T})$ as a diagonal subalgebra. Both corresponding equivalence relations are equivalent: they have the same quotient space $\bf T$. The first one is given by the trivial covering map from ${\bf T}\times \{1,\ldots,n\}$ onto ${\bf T}$ while the second is given by the covering map $z\to z^n$ from $\bf T$ onto $\bf T$.

2. Let $\varphi: G\rightarrow H$ be a continuous homomorphism of locally compact abelian groups $G,H$. Then $G$ acts continuously on $H$ by left multiplication and we can form the crossed product C*-algebra $G\,\lsd C_0(H)$. By dualizing, we get $\hat\varphi: \hat H\rightarrow \hat G$ and the crossed product C*-algebra $\hat H\,\lsd C_0(\hat G)$. The Fourier transform gives an isomorphism of these C*-algebras.

If $G$ is discrete and $\varphi$ is one-to-one, $C_0(H)$ is a diagonal subalgebra. Similarly, if $\hat H$ is discrete and $\hat\varphi$ is one-to-one, $C_0(\hat G)$ is another diagonal subalgebra. Both conditions happen simultaneously if $G$ is discrete, $H$ is compact, $\varphi$ is one-to-one and has dense range. There are such examples where $C_0(H)$ and $C_0(\hat G)$ are not isomorphic.
\begin{example}
$G={\bf Z}^2, H={\bf R}/{\bf Z}, \varphi(m,n)=\alpha m+\beta n+{\bf Z}$
where $(1,\alpha,\beta)$ are linearly independent over $\bf Q$.
\end{example}

\subsection{Variations on the cross} The cross consists of the graph of the functions $y=x$ and $y=-x$ on the domain $[-1,1]$. There are several ways to deal with the singular point $(0,0)$. We shall present three of them.

The C*-algebra $A=C([0,1])\otimes M_2({\bf C})$ has the obvious diagonal
$B=C([0,1])\otimes D_2$.
It is instructive to look at the pairs $(A_i, B_i=A_i\cap B)$, where $A_i$ is one of the following subalgebras of $A$.
\vskip 3mm
$\begin{array}{cc}
 A_1=&\{f\in A:f(0)=\left(\begin{array}{cc}a & b \\b & a\end{array}\right)\}\\
 &\\
 A_2=&\{f\in A:f(0)=\left(\begin{array}{cc}a & a \\a & a\end{array}\right)\}\\
 &\\
 A_3=&\{f\in A:f(0)=\left(\begin{array}{cc}a & 0 \\0 & a\end{array}\right)\}\\
 &\\
 A_4=&\{f\in A:f(0)=\left(\begin{array}{cc}a & 0 \\0 & b\end{array}\right)\}
\end{array}$  

\vskip 3mm
\subsubsection{A groupoid of germs} In the first example, $B_1$ is a Cartan subalgebra of $A_1$ which does not have the unique extension property. Indeed, the  states  $f\mapsto a\pm b$  both extend the pure state 
$f\mapsto a$  of $B_1$.
 The pair $(A_1,B_1)$ can be realized as $(C^*(G), C([-1,1]))$, where $G$ is the groupoid of the action of the group ${\bf Z}/2{\bf Z}$ on $[-1,1]$ by the map $Tx=-x$. This groupoid can also be described as the groupoid of germs of the pseudogroup generated by $T$:
$$G=\{(\pm x,\pm 1, x), x\in [-1,1]\}.$$
It is topologically principal but not principal: the isotropy is trivial at $x\not=0$ and the isotropy group at $x=0$ is ${\bf Z}/2{\bf Z}$. It is an elementary example of an orbifold. The C$^*$-algebra $A_1$  is a CCR algebra with non-Hausdorff spectrum. I owe to A. Kumjian the observation that $A_1$ does have a diagonal. Indeed, it is isomorphic to $A_4$, which has $B_4$ as a diagonal subalgebra.

\subsubsection{A branched covering}

The subalgebra $B_2$ is maximal abelian in $A_2$. However, it  does not satisfy the conditions $(1)$ and $(4)$ of the definition of a Cartan subalgebra. Its main defect is to be contained in the ideal $f(0)=0$. The $C^*$-algebra $A_2$ can be realized as the C*-algebra of a non-\'etale principal groupoid, namely the equivalence relation $R$ associated to the previous groupoid $G$. Endowed with the product topology of $[-1,1]\times [-1,1]$, it is a proper groupoid. It has the Haar system: $$\int f d\lambda^x=f(x,x)+f(x,-x).$$
The quotient map, which can be realized as the map $x\to |x|$ from $[-1,1]$ to $[0,1]$ is an elementary example of a branched covering (see \cite{fox:coverings}). The construction given here appears in \cite{dea:thesis}.  The C*-algebra $A_2=C^*(R,\lambda)$ does not contain $C([-1,1])$ as a subalgebra because the diagonal of $R$ is not open; however, it contains $B_2=C([-1,1]\setminus\{0\})$. The quotient map $G\rightarrow R$ gives the inclusion $A_2\subset A_1$. Note also that $A_1$ is the unitization of $A_2$. It can be readily checked or deduced from \cite{mrw:continuous} that $A_2$ is a continuous trace C*-algebra.

\subsubsection{An \'etale equivalence relation}

The subalgebra $B_3$ is a diagonal subalgebra in $A_3$. It is realized by the same equivalence relation $R$ as above, but endowed with a finer topology which makes it \'etale. Following M.~Molberg \cite{mol:AF}, we consider the topology $\tau$ generated by the product topology and the diagonal $\{(x,x), x\in [0,1]\}$.  Then $R_\tau$ is \'etale but no longer proper.

\section{Orbifolds and non-Hausdorff manifolds}

Moerdijk and Pronk have introduced \cite{mp:orbi} the notion of an orbifold groupoid (here, we only need the topological structure of the spaces, not their differential structure):

\begin{definition} An {\it orbifold groupoid} is  a proper, effective, \'etale, second countable, locally compact and Hausdorff groupoid.
\end{definition}
 
 They arise in the following related situations:
\begin{itemize} 
\item orbifolds;
\item foliated manifolds for which all the leaves are compact with finite holonomy.
\end{itemize}
The explicit constructions involve some choices but provide equivalent groupoids (see \cite{mm:intro}). It seems appropriate to define a (topological) orbifold as an equivalence class of proper (effective) groupoids. The example 2.2.1 of the previous section is an elementary example of an orbifold groupoid. The general case keeps some of the features of this elementary example: if $G$ is an orbifold groupoid, its $C^*$-algebra is a CCR algebra which admits $C_0(G^{(0)})$ as a Cartan subalgebra. There is an intriguing  link between the orbifold groupoid $G=[-1,1]\rsd\,{\bf Z}/2{\bf Z}$ of example 2.2.1  and the non-Hausdorff manifold obtained as the quotient of $[0,1]\times \{0,1\}$ by the equivalence relation which identifies the two copies of $(0,1]$: 
$$R=\{(z,(i,j))\in [0,1]\times(\{0,1\}\times \{0,1\}): i=j\quad\hbox{if}\quad z=0\}.$$
This \'etale equivalence relation gives the diagonal $B_4\subset A_4$ and we have seen that $A_1$ and $A_4$ are isomorphic C$^*$-algebras. The isotropy has been eliminated at the expense of non-Hausdorffness. 
\vskip 3mm
Non-Hausdorff manifolds are a rich source of examples of CCR algebras. Let us define a locally compact space as topological space $Y$ such that every point has a compact Hausdorff neighborhood. Such a space is $T_1$. Let us define a desingularization of $Y$ as a surjective local homeomorphism $\pi:X\rightarrow Y$, where $X$ is a Hausdorff locally compact space. Then, the graph $R$ of the equivalence relation $\pi(x)=\pi(x')$ on $X$, endowed with the product topology of $X\times X$ is an \'etale equivalence relation. It is proper if and only if $Y$ is Hausdorff. The C$^*$-algebra $C^*(R)$ is CCR and its spectrum is homeomorphic to $Y$ (see for example \cite{orl:principal} for these facts). This provides a convenient way to construct CCR algebras with arbitrary locally compact spectrum. Let $Y$ be a topological space which is $T_1$. One says that two points of $Y$ are separated if they have disjoint neighborhoods and that a point $y\in Y$ is Hausdorff if it is separated from any other point. The set of Hausdorff points of the spectrum of a separable CCR algebra is a dense $G_\delta$. In many cases (for example, when the spectrum is compact (not necessarily Hausdorff) or in the case of the C$^*$-algebra of a connected and simply connected nilpotent Lie group), the interior of this set is dense. However, in \cite{dix:separes}, Dixmier gives an example of a separable CCR algebra  such that the set of Hausdorff points of its spectrum has an empty interior. Here is an easy construction of a similar algebra inspired by Dixmier's example (I do not know whether the algebras are the same) and also by the example $A_4$. Let $Z$ be a Hausdorff locally compact space and let $\{z_n, n\in {\bf N}\}$ be a countable subset of $Z$. Then
$$R=\{(z,(i,j))\in Z\times({\bf N}\times{\bf N}): i=j\quad\hbox{if}\quad z\in\{z_i,z_j\}\}$$
is an open subgroupoid of $ Z\times({\bf N}\times{\bf N})$. Therefore, it is an \'etale equivalence relation over $X=Z\times{\bf N}$. The quotient space $Y=X/R$ can be described as the disjoint union of $Z$ and ${\bf N}$. The quotient map sends $(z,i)\in X$ to $z$ if $z\not=z_i$ and to $i$ if $z=z_i$. Open subsets of $Z$, where a finite number of $z_i$'s have been replaced by $i$ form a base for the quotient topology. The space $Y$ is locally compact. If none of the $z_i$'s are isolated, the set of Hausdorff points is
$Z\setminus\{z_i, i\in {\bf N}\}$. It has an empty interior if $\{z_i, i\in {\bf N}\}$ is dense in $Z$. As said before, the C$^*$-algebra $C^*(R)$ is CCR and has $Y$ as its spectrum.

\section{The non-Hausdorff case} According to Theorem 1.3, Cartan subalgebras in C$^*$-algebras are intimately related to topologically principal \'etale groupoids. These groupoids arise as groupoids of germs. However, the Hausdorffness condition required in the theorem is a severe restriction. It is still possible to define the reduced C$^*$-algebra $C^*_r(G)$ (and $C^*_r(G,E)$) when $G$ is an \'etale locally compact non-Hausdorff groupoid. By definition, an element of $C_c(G)$ is a function of the form $f=\sum_{i=1}^n \tilde f_i$, where $f_i\in C_c(U_i)$ for some open Hausdorff subset $U_i\subset G$ and $\tilde f_i$ is its extension by zero to $G$. Then, the definitions are just as above. A function $f$ in $C_c(G)$ is not necessarily continuous on $G$. In particular, its restriction to $G^{(0)}$ is not necessarily continuous. Thus, the existence of a conditional expectation onto $C_0(G^{(0)})$ is problematic. The subalgebra $C_0(G^{(0)})$ may also fail to be maximal abelian.
Shortly after the workshop, R. Exel gave me an example of a non-Hausdorff groupoid of germs $G$ such that the subalgebra $C_0(G^{(0)})$ is not maximal abelian (see \cite{exe:example}). I then realized that this example is related to an earlier example of G. Skandalis which appears in \cite{ren:ideal} and which I reproduce below. Skandalis' purpose is different: it shows that the C$^*$-algebra of a minimal foliation is not necessarily simple when the holonomy groupoid is non-Hausdorff. However, both pathologies are based on the same fact. 

Let $g_1, g_2$ be homeomorphisms of the circle $\bf T$ having for fixed points set respectively the oriented arcs $[a,b]$ and $[b,a]$, where $a,b$ are distinct points of $\bf T$. Since these homeomorphisms commute, they define an action of ${\bf Z}^2$ on $\bf T$ such that $(m,n)$ acts as $g_1^mg_2^n$. Let $[g(x),g,x]$ denote the germ of a homeomorphism $g$ at $x\in{\bf T}$. Let $G$ be the groupoid of germs of the $g_1^mg_2^n$'s. The only points which have non trivial isotropy are $a$ and $b$. The isotropy subgroups at $a$ and $b$ are isomorphic to ${\bf Z}^2$.  By construction, $S(m,n)=\{[g_1^mg_2^n(x),g_1^mg_2^n,x]: x\in{\bf T}\}$ is an open bisection of $G$.  The function 
$$f={\bf 1}_{S(0,0)}-{\bf 1}_{S(1,0)}-{\bf 1}_{S(0,1)}+{\bf 1}_{S(1,1)}$$
belongs to $C_c(G)$, hence to $C^*_r(G)$. It vanishes outside the finite set 
$$\{[a,(m,n),a], [b,(m,n),b], m,n=0,1\}$$ 
but takes the values $\pm 1$ on this set. Since this function has its support contained in the isotropy group bundle $G'$, it commutes with every element of $C({\bf T})$, however it does not belong to $C({\bf T})$.
\vskip 5mm
{\it Acknowledgements.} I thank the participants of the workshop for stimulating discussions and comments, in particular R.~Archbold and R.~Exel. I also benefitted greatly from the help of A.~Kumjian.

\bibliographystyle{amsalpha}

\end{document}